\documentclass[12pt]{amsart}       %{proc-l}                %[12pt]{amsart}
\usepackage{amsmath, amstext, amsbsy}
\newtheorem{theorem}{Theorem}[section]
\newtheorem{lemma}[theorem]{Lemma}
\theoremstyle{definition}

\theoremstyle{remark}

\numberwithin{equation}{section}

\renewcommand{\a}{\alpha}

\newcommand{\ep}{\epsilon}

\begin{document}
%\magstep2
%Topmatter

%One author
\title[Quantum Kac-Moody algebras]
    {Quantum Kac-Moody algebras
    \\and vertex representations %<title line 2>
    }
\author{Naihuan Jing}
\address{Department of Mathematics,
   North Carolina State Univer\-sity,
   Ra\-leigh, NC 27695-8205,
   U.S.A.}
\email{jing@eos.ncsu.edu}
\thanks{Research supported in part by NSA grant MDA904-97-1-0062.}
%End One author

%\date{}

%End topmatter

\begin{abstract}
We introduce an affinization of the quantum Kac-Moody algebra associated
to a symmetric generalized Cartan matrix.
Based on the affinization, we construct a representation of the
quantum Kac-Moody algebra
by vertex operators from bosonic fields. We also obtain a combinatorial
indentity about Hall-Little\-wood polynomials.
\end{abstract}

\maketitle

\centerline{To appear in Lett. Math. Phys.}

\section{Introduction} \label{S:intro}

Vertex representations of quantum affine algebras were first studied in
\cite{FJ}\cite{J1} for irreducible level one
representations of quantum affine
algebras of both untwisted and twisted types. With the latest realizations of
quantum affine algebras $U_q(G_2^{(1)})$ \cite{J3} and
$U_q(F_4^{(1)})$ \cite{J4} all quantum affine
algebras have been constructed explicitly at least for level one
and some for other levels. An overview of the constructions in terms of
quantum $Z$-algebras was given in \cite{J4} as well as the
historic developments circled around various constructions.

In an abstract setting, Lusztig established the theory of
highest weight modules for arbitrary symmetrizable quantum Kac-Moody
algebras \cite{L}. These are done by prescribing the actions of the Serre
generators of the quantum Kac-Moody algebras \cite{D, Jb}
paralleled to the classical cases
of $q=1$. On the other hand in the classical case of $q=1$, I. Frenkel \cite{F}
had constructed a vertex representation for the Kac-Moody Lie
algebra associated to a symmetric generalized Cartan matrix
in order to study root multiplicities of the
hyperbolic Kac-Moody algebra.

In this paper we will give a quantum generalization of I. Frenkel's
construction
for quantum Kac-Moody algebras. We first imbed the quantum
Kac-Moody algebra into an affinization of
the quantum Kac-Moody algebra. Then we construct a (level one)
vertex representation
of the quantum affinization.
Our construction gives a nontrivial
 combinatorial identity (cf. \ref{Ser2}), which is interesting by itself.
This identity was trivial in the classical case and
was avoided by using usual vertex operator techniques in \cite{F}.
We proved the quantum Serre relations by a spiral method using the
representation.

The paper is organized as follows. In section two we review the basic
definition of quantum enveloping algebras associated to generalized Cartan
matrices. In section three we introduce a quantum Heisenberg algebra
associated to the root lattice of the Kac-Moody algebra and the representation
is stated thereafter. In sections four and five we set out to prove that the
construction gives a realization of the quantum Kac-Moody algebra. The proof
of the Serre relations is obtained by showing a stronger identity
of vertex operators, and this latter identity is equivalent to an identity
of symmetric functions. Finally we remark about the connection of the
combinatorial identity (\ref{Ser2}) with
Hall-Littlewood polynomials.

\section{Quantum Kac-Moody algebras and their affinizations} \label{S:alg}

The notation of Kac-Moody algebras essentially follows \cite{K}.
Let $g(A)$ be the Kac-Moody Lie algebra associated to the generalized
Cartan matrix $A=(a_{ij})$, where $a_{ii}=2, a_{ij}=a_{ji}\leq 0$ for $i, j=1,
\cdots, l$. Let $Q=\mathbb Z\a_1\oplus\cdots
\oplus \mathbb Z\a_l$ be the root lattice
generated by the simple roots $\a_i$, $Q^\vee=\mathbb Z h_1\oplus\cdots\oplus
\mathbb Z h_l$ the coroot lattice generated by the coroots
$h_i$, and $P^\vee=\mathbb Z\lambda_1\oplus\cdots\oplus
\mathbb Z\lambda_l$ the weight lattice generated by the fundamental weights
$\lambda_i$.
The normalized invariant bilinear form $(\ \ |\ \ )$ on $g(A)$
satisfies the following properties:
$$(\a_i|\a_i)=2, \ \ (\a_i|\a_j)=a_{ij}. $$

Let $q$ be a generic complex number (that is, $q$ is not a root of unity).
The quantum Kac-Moody algebra is an associative algebra $U_q=U_q(g(A))$
over ${\mathbb C}(q^{1/2})$ generated by elements
$e_i, f_i, q^{h}$, $i=1, \cdots, l$, $h\in P^*=Hom(P, \mathbb Z)$ with the
relations:
\begin{align*}
   q^h&=1 \quad \hbox{for $h=0$,} \\
   q^{h+h'}&=q^h q^{h'} \quad \hbox{for $h,h'\in {P}^*$,} \\
   q^h e_i q^{-h} &= q^{(h | \alpha_i)} e_i
      \quad \hbox{and} \quad
      q^h f_i q^{-h} = q^{-(h | \alpha_i )}f_i, \\
   [e_i,f_j]&=\delta_{ij}\frac{t_i-t_i^{-1}}{q-q^{-1}}, \\
   \sum_{k=0}^{1-a_{ij}}
      \frac{(-1)^k}{[k]![b-k]!}& e_i^{k} e_j e_i^{b-k}=0,
      \quad
      \sum_{k=0}^{b}
      \frac{(-1)^k}{[k]![b-k]!} f_i^{k} f_j f_i^{b-k}=0
      \quad \mbox{for} \; i \neq j. \\
\end{align*}
where $t_i=q^{h_i}$ and
\begin{equation*}
      [k]=\frac{q^k-q^{-k}}{q-q^{-1}}, \quad
      [k]!=[1][2]\cdots [k], \quad
      \; \hbox{ and } \;
      t_i=q^{h_i}%=q^{\a_i}.
\end{equation*}
   The associative algebra $U_q$ has a Hopf algebra structure
   with the following comultiplication.
\begin{equation*}
 \Delta( q^h)=q^h\otimes q^h, \;
      \Delta(e_i)=e_i \otimes 1 + t_i \otimes e_i, \;
      \Delta(e_i)=f_i \otimes t_i^{-1} + 1 \otimes f_i.
\end{equation*}

Generalizing the Drinfeld loop algebra realization of quantum affine algebras,
we introduce the following associate algebra, which is a $q$-deformation of the
affinization of Kac-Moody algebra.

The $q$-affinization algebra is an associative algebra
$U_q(\widehat{g(A)})$ generated by
elements
$x^{\pm}_{in}$, $a_{ik}$, $K_i^{\pm}$, $q^d$ and $\gamma$ for
 $n\in \mathbb Z, k\in
\mathbb Z-{0}$ with the following defining relations:
\begin{align} \label{E:R1}
 [\gamma^{\pm 1/2}, u]&=0 \ \ \mbox{for all} \ u\in \textstyle {\bf U},\\
\mbox{} [a_{ik}, a_{jl}]&=\delta_{k+l,0}
\displaystyle\frac {[(\a_i|\a_j)k]}{k}
\displaystyle\frac {\gamma^k-\gamma^{-k}}{q-q^{-1}}, \label{E:R2}\\
\mbox{} [a_{ik}, K_j^{\pm 1}]&=[q^{\pm d}, K_j^{\pm 1}]=0, \label{E:R3}\\
 q^d x_{ik}^{\pm} q^{-d}&=q^k x_{ik}^{\pm }, \ \
q^d a_{il} q^{-d}=q^l a_{il}, \label{E:R4}\\
 K_i x_{jk}^{\pm} K_i^{-1}&=q^{\pm (\alpha_i|\alpha_j)} x_{jk} ^{\pm},
\label{E:R5}\\
\mbox{} [a_{ik}, x_{jl}^{\pm}]&=\pm \displaystyle\frac {[(\a_i|\a_j)k]}{k}
\gamma^{\mp |k|/2} x_{j,k+l}^{\pm},            \label{E:R6}
\end{align}
\begin{align}
\prod_{r=0}^{-(\a_i|\a_j)+1}
(z-&q^{\pm(\a_i|\a_j)+2r}w)x_i^{\pm}(z)x_j^{\pm}(w) \nonumber \\
          +&\prod_{r=0}^{-(\a_i|\a_j)+1}(w-q^{\pm(\a_i|\a_j)+2r}z)
x_j^{\pm}(w)x_i^{\pm}(z) =0, \label{E:R7}
\\
\mbox{}
\ [x_i^+(z),x_j^-(w)]
             =\frac{\delta_{ij}}{q-q^{-1}}&
                 \left(
                 \psi_i(w\gamma^{1/2})
                 \delta(\frac{w\gamma}{z})
                 -\varphi_i(w\gamma^{-1/2})
                  \delta(\frac{w\gamma^{-1}}z)
                 \right)                         \label{E:R8}
\end{align}
where
$x_i^{\pm}(z)=\sum_{n\in {\bf Z}}x_{i,n}z^{-n-1}$, $\psi_{im}$ and
$\varphi_{im}$ $(m\in {\bf Z}_{\ge 0})$
are defined by
\begin{align}
 \sum_{m=0}^{\infty} \psi_{im} z^{-m}
 &=K_i \textstyle {exp} \left( (q-q^{-1}) \sum_{k=1}^{\infty} a_{ik}
z^{-k}
\right), \\
 \sum_{m=0}^{\infty} \varphi_{i,-m} z^{m}
 &=K_i^{-1} \textstyle {exp} \left(- (q-q^{-1}) \sum_{k=1}^{\infty}
a_{i,-k} z^{k}
\right),
\\
 \sum_{r=0, \sigma\in S_m}^{m=1-a_{ij}}(-1)^r &
\left[\begin{array}{c} m\\r\end{array}\right]\sigma .x_i^{\pm}(z_1)\cdots
x^{\pm}_{i}(z_r)  \nonumber\\
&\qquad\qquad \cdot x^{\pm}_{j}(w) x^{\pm}_{i}(z_{r+1})\cdots x^{\pm}_{i}
(z_m)=0 . \label{E:R9}
\end{align}
where the symmetric group $S_m$ acts on  $z_i$'s
by permuting their indices.

\begin{theorem} \label{T:1}
The mapping $e_i\mapsto x^+_{i, 0}, f_i\mapsto x_{i, 0}^-,
t_i\mapsto K_i$ defines an imbedding of the quantum Kac-Moody algebra
$U_q(g(A))$ into the quantum affinization algebra $U_q(\widehat{g(A)})$.
\end{theorem}
\begin{proof} It is easy to see that $X_{i0}^{\pm}$ and $K_i$
 satisfy the relations
of $e_i, f_i, t_i$. The only thing needs explanation is the Serre relation.
This follows from (\ref{E:R9}) by taking the constant coefficient of
$w, z_1, \cdots, z_m$ and noting the invariance of each term under
$S_m$.
\end{proof}

\section{Fock space representations}   \label{S:fock}

The quantum Heisenberg algebra is the associative algebra generated by $a_i(n)$
and the central element $\gamma^{\pm 1}$ satisfying the relations
(\ref{E:R2}).We
will consider the level one $q$-Heisenberg algebra, that is $\gamma=q$.
Let $a_i(m)$ ($i=1, \cdots, l$, $m\in \mathbb Z$) be the operators satisfying:
\begin{equation}
 [a_i(m), a_j(n)]=\delta_{m+n, 0}\frac{[(\a_i|\a_j)m]}{m}[m].
\end{equation}
where we take the limit of $m\mapsto 0$ if the mode $m=0$.

   We define the Fock module ${\mathcal F}$ as the tensor product of the
   symmetric
   algebra generated by $a_i(-n)$ ($n\in \mathbb N$) and
   the twisted group
   algebra $\mathbb C\{Q\}$ generated by $e^{\a}, \a\in Q$ subject to
relation:
$$e^{\alpha_1}e^{\alpha_2}=(-1)^{(\a_i|\a_j)}e^{\alpha_2}e^{\alpha_1}.
$$

Let the element
$1$ be  the vacuum state with the actions of $\a_i(n)$:
   $$ a_i(n).1 = 0 \quad (n>0).
  $$
The element $a_i(0)$ act as differential operators  by
$$a_i(0)e^{\alpha}=(\alpha_i|\alpha)e^{\alpha}.$$
Clearly this defines a representation of the $q$-Heisenberg algebra for
$\gamma=q$.

  As usual we define the normal product as the ordered product by moving
annihilation operators $a_i(n), a_i(0), a_i(0)$ $(n>0$) to the left.

   Let us introduce the following vertex operators.
   \begin{equation}
   \begin{aligned}
   X_i^{\pm}(z)&=
         \exp ( \pm \sum^{\infty}_{n=1}
                 \frac{a_i(-n)}{[n]} q^{\mp \frac{n}{2}} z^n)
         \exp ( \mp \sum^{\infty}_{n=1}
                 \frac{a_i(n)}{[n]} q^{\mp \frac{n}{2}} z^{-n})
         e^{\a_i} z^{\mp a_i(0)}\\
          &=\sum_{n\in \mathbb Z}X_i^{\pm}(n)z^{-n-1}.
   \end{aligned}
   \end{equation}
 %  \bigskip
   \noindent
\begin{theorem} \label{T:2} The space
    ${\mathcal F}$ is a $U_q(\widehat{g(A)})$-module of level one
    under the action defined by
    $ \gamma \longmapsto q,
       K_i \longmapsto q^{a_i(0)}$,
    $ a_{im} \longmapsto a_i(m),
       q^d \longmapsto q^{{d}}$, and
\begin{equation*}
x_{i, n}^{\pm} \longmapsto X_{i}^{\pm}(n).
\end{equation*}
\end{theorem}

The proof of the theorem will occupy this and next two sections.
We start by recalling the $q$-binomial functions $(1-x)^n_{q^2}$ introduced in
\cite{J4}:
$$
(1-x)^a_{q^2}=\frac{(q^{-a+1}x;q^2)_{\infty}}{(q^{a+1}x;q^2)_{\infty}}.
$$
In particular for $a\in N$ we have:
$$
(1-x)^a_{q^2}=(1-q^{-a+1}x)(1-q^{-a+3}x)\cdots (1-q^{a-1}x)=(q^{-a+1}x;q^2)_a.
$$
For convenience we will also use the following $q$-analogue of binomial
functions:
\begin{equation}\label{E:N1}
(z-w)^a_{q^2}=(1-\frac wz)_{q^{2}}^{a}z^a.
\end{equation}

Under our symmetric $q$-numbers, the $q$-binomial theorem becomes:
\begin{equation}\label{binom1}
\sum_{r=0}^{n}\begin{bmatrix} n\\r\end{bmatrix}(-z)^{r}=(q^{1-n}z;q^2)_n.
\end{equation}
In particular we have
\begin{equation}\label{binom2}
\sum_{r=0}^{n}\begin{bmatrix} n\\r\end{bmatrix}(-q^i)^{r}=0,
\qquad\mbox{for } i=1-n, 3-n, \cdots, n-1.
\end{equation}

We now prove the theorem by checking that the action satisfies the
relations of the $q$-affinization algebra. It is clear that the relations
(\ref{E:R1}-\ref{E:R5}) are true by the construction. Relation
(\ref{E:R6}) follows from the definition of $X^{\pm}_i(z)$.
So we only need to show the
relations (\ref{E:R7}-\ref{E:R8}) and the Serre relation (\ref{E:R9}).

The operator product expansions (OPE)
for $X_i^{\pm}(z)$ are computed as follows (cf.
\cite{J4}):
\begin{equation}\label{OPE1}
\begin{aligned}
X_i^{\pm}(z)X_j^{\pm}(w)&=:X_i^{\pm}(z)X_j^{\pm}(w):(1-q^{\mp 1}\frac
wz)_{q^2}^{(\a_i|\a_j)}z^{(\a_i|\a_j)},\\
X_i^{\pm}(z)X_j^{\mp}(w)&=:X_i^{\pm}(z)X_j^{\mp}(w):(1-\frac
wz)_{q^2}^{-(\a_i|\a_j)}z^{-(\a_i|\a_j)}.
\end{aligned}
\end{equation}
where the contraction functions like $(1-q^{\mp 1}\frac
wz)_{q^2}^{(\a_i|\a_j)}z^{(\a_i|\a_j)}$ are understood as power
series in $w/z$. Using the notation (\ref{E:N1}) we can combine the OPE
into a self-suggested form:
\begin{align}\label{OPE2}
X_i^{\epsilon}(z)X_j^{\epsilon'}(w)&=:X_i^{\epsilon}(z)X_j^{\epsilon'}(w):(z-
q^{(\ep+\ep')/2}
w)_{q^2}^{(\a_i|\a_j)},\\
:X_i^{\epsilon}(z)X_j^{\epsilon'}(w):&=(-1)^{(\a_i|\a_j)}
:X_j^{\epsilon}(z)X_i^{\epsilon'}(w):. \label{OPE3}
\end{align}

The operator product expansions imply immediately the commutation relations
(\ref{E:R7}). We will find out that there exist more detailed commutation
relations
between $X_i^{\pm}(z)$ and $X_j^{\pm}(w)$.

The case of $i=j$ in the relations (\ref{E:R8}) follows from the case of
$U_q(\hat{sl}_2)$ \cite{J2}. When $i\neq j$ we have
\begin{align*}
[X_i^{\pm}(z), X_j^{\mp}(w)]&=:X_i^{\pm}(z)X_j^{\mp}(w):
\left((1-\frac
wz)_{q^2}^{(\a_i|\a_j)}z^{-(\a_i|\a_j)}\right.\\
&\qquad\qquad \left. (-1)^{(\a_i|\a_j)}(1-\frac
zw)_{q^2}^{(\a_i|\a_j)}w^{-(\a_i|\a_j)}\right)\\
&=0
\end{align*}
since the contraction functions $(z-w)_{q^2}^{(\a_i|\a_j)}=
(-1)^{(\a_i|\a_j)}(w-z)_{q^2}^{(\a_i|\a_j)}$ as $q$-polynomials in $z, w$
(cf. \ref{E:N1}).

\section{q-difference operators and commutation relations}

The commutation relations of operators $X_i^{\pm}(z)$ and $X_j^{\pm}(w)$
can be improved by using $q$-difference operators.

Let $f(x)$ be a function. The $q$-difference operator $D_q(x)$ is defined
by
$$D_q(x;q)f(x)=\frac{f(xq)-f(xq^{-1})}{(q-q^{-1})x}.
$$
We will simply write $D_q$ if the variable is clear from the context.

\begin{lemma} \label{diff1}
$D_q(1-x)_{q^2}^n=-[n](1-x)_{q^2}^{n-1}.$
\end{lemma}

\begin{proof} We compute directly that
\begin{align*}
D_q(1-x)_{q^2}^n&=\frac 1{(q-q^{-1})x}((q^{-n+2}x;q^2)_n-(q^{-n}x;q^2)_n)\\
&=\frac 1{(q-q^{-1})x}(q^{-n+2}x; q^2)_{n-1}(1-q^nx-1+q^{-n}x)\\
&=-[n](q^{-n+2}x; q^2)_{n-1}
\end{align*}
\end{proof}

The following result will be useful for what follows.

\begin{lemma} \label{diff2}
$$D_q^nf(x)=\frac 1{(q-q^{-1})^nx^n}\sum_{i=0}^n(-1)^iq^{(n-1)(-n+2i)/2}
\begin{bmatrix}n\\i\end{bmatrix}f(q^{n-2i}x).
$$
\end{lemma}
\begin{proof} By induction on $n$ we have
\begin{align*}
&D_q^nf(x)=\frac{D_{q}^{n-1}f(qx)-D_{q}^{n-1}f(q^{-1}x)}{(q-q^{-1})x}\\
 &\frac
1{(q-q^{-1})^nx^n}\left(q^{-n+1}\sum_{i=0}^{n-1}(-1)^iq^{(n-2)(-n+1+2i)/2}
\begin{bmatrix}n-1\\i\end{bmatrix}f(q^{n-2i}x)\right.\\
&\qquad \left. -q^{n-1}\sum_{i=0}^{n-1}(-1)^iq^{(n-2)(-n+1+2i)/2}
\begin{bmatrix}n-1\\i\end{bmatrix}f(q^{n-2-2i}x)\right)\\
&=\frac 1{(q-q^{-1})^nx^n}\left(q^{-\binom{n}{2}}f(q^nx)+
(-1)^nq^{\binom{n}{2}}f(q^{-n}x)\right.\\
&\qquad \left.+\sum_{i=1}^{n-1}(-1)^iq^{-\binom{n}{2}+i(n-1)}
\left(q^{-i}
\begin{bmatrix}n-1\\i\end{bmatrix}+
q^{n-i}
\begin{bmatrix}n-1\\i-1\end{bmatrix}\right)f(q^{n-2i}x)\right)\\
&=\frac 1{(q-q^{-1})^nx^n}\sum_{i=0}^n(-1)^iq^{-\binom{n}{2}+i(n-1)}
\begin{bmatrix}n\\i\end{bmatrix}f(q^{n-2i}x).
\end{align*}
\end{proof}

To state the improved commutation relation we define:
$$
X^{\pm}_{\a_i+\a_j, r}(z)=:X_i^{\pm}(z)X_j^{\pm}(zq^{(\a_i|\a_j)+2r+1}):
z^{(\a_i|\a_j)+1},
$$
for $r=0, \cdots, -(\a_i|\a_j)-1$.
\begin{theorem} \label{T:3} Let $m=1-(\a_i|\a_j)$.
On the Fock space $V$ we have the following commutation
relation:
\begin{align*}
&X_i^{\pm}(z)X_j^{\pm}(w)-X_j^{\pm}(w)X_i^{\pm}(z)=0, \qquad\mbox{if
$(\a_i|\a_j)=0$};\\
&(z-q^{\pm(\a_i|\a_j)}w)X_i^{\pm}(z)X_j^{\pm}(w)+
(w-q^{\pm(\a_i|\a_j)}z)X_j^{\pm}(w)X_i^{\pm}(z)=0, \\
&\hskip 2.5in \qquad\mbox{if
$(\a_i|\a_j)=-1$};\\
&(z-q^{\pm(\a_i|\a_j)}w)X_i^{\pm}(z)X_j^{\pm}(w)+
(w-q^{\pm(\a_i|\a_j)}z)X_j^{\pm}(w)X_i^{\pm}(z)\\
&=:X_i^{\pm}(z)X_j^{\pm}(w):\frac{1}{(q-q^{-1})\cdots(q^{m-2}-q^{-m+2})}\\
&\qquad\qquad \sum_{r=0}^{m-2} q^{-\binom{m-2}{2}+r(m-3)}\begin{bmatrix} m-2
\\r\end{bmatrix}
\delta(wq^{m-2-2r}/z), \mbox{if $(\a_i|\a_j)\leq -2$}.
\end{align*}
\end{theorem}
\begin{proof} We only need to show the case of $+$.
For $m=1-(\a_i|\a_j)\geq 2$
it follows from (\ref{OPE1}-\ref{OPE2}) that
\begin{align*}
&(z-q^{(\a_i|\a_j)}w)X_i^{+}(z)X_j^{+}(w)+
(w-q^{(\a_i|\a_j)}z)X_j^{+}(w)X_i^{+}(z)\\
&=:X_i^{+}(z)X_j^{+}(w):\left(\frac 1{(z-q^{3-m}w)(z-q^{5-m}w)\cdots
(z-q^{m-3}w)}\right.\\
&\qquad \left. +\frac{(-1)^{m-2}}{(w-q^{3-m}z)(w-q^{5-m}z)\cdots
(w-q^{m-3}z)}\right)\\
&=:X_i^{+}(z)X_j^{+}(w):\frac1{[m-2]!z^{m-2}}\left(D_q^{m-2}(w/z)(
\frac 1{1-w/z}
+\frac{z/w}{1-z/w})\right)\\
 &=\frac{:X_i^{+}(z)X_j^{+}(w):}{(q-q^{-1})\cdots (q^{m-2}-q^{2-m})z^{m-2}}
\sum_{r=0}^{m-2}(-1)^rq^{(m-3)(-m+2+2r)/2}\\
&\qquad\qquad\qquad \cdot\begin{bmatrix}m-2\\r\end{bmatrix}\delta
(q^{m-2-2r}w/z)\\
 &=\frac 1{(q-q^{-1})^{m-2}[m-2]}\sum_{r=0}^{m-2}q^{-\binom{m-2}{2}+
r(m-3)}
X_{\a_i+\a_j, r}^+(z)\delta (q^{m-2-2r}\frac wz).
\end{align*}
\end{proof}

It is clear that the commutation relation (\ref{E:R7}) is a consequence of
Theorem (\ref{T:3}).

\section{Quantum Serre relations}

The proof of Serre relations shows the power of our vertex representation.
As an application we obtain a non-trivial combinatorial identity.

\begin{lemma}\label{T:4} Suppose that there exists a root $\a_k$
such that $(\a_k|\a_i)=0$ and $(\a_k|\a_j)\neq 0$, that is, $\a_k$ joins
with $\a_j$ but not with $\a_i$. Then the $q$-affinized
Serre relation (\ref{E:R9}) associated
with $a_{ij}$ is equivalent to the following identity:
\begin{equation}\label{E:R10}
 \sum^{m=1-a_{ij}}_{r=0, \sigma\in S_m}(-1)^r
\left[\begin{array}{c} m\\r\end{array}\right]\sigma.x_{i}^{\pm}(z_1)\cdots
x^{\pm}_{i}(z_r)x^{\pm}_{j, l}x^{\pm}_{i}(z_{r+1})\cdots x^{\pm}_{i}(z_m)=0
\end{equation}
for some $l\in \mathbb Z$ and $\sigma$ acts on indices of $z_i$'s.
\end{lemma}
\begin{proof} The relation (\ref{E:R10}) is a special case of (\ref{E:R9}),
so we only need to show the converse direction. Assume that (\ref{E:R10})
holds for some $l\in \mathbb Z$.
Taking commutator with $a_k(s)$ to (\ref{E:R10})
using (\ref{E:R6}),
we see immediately by induction
that
\begin{equation}\label{E:R11}
\sum^{m=1-a_{ij}}_{r=0, \sigma\in S_m}(-1)^r
\left[\begin{array}{c} m\\r\end{array}\right]\sigma.x_{i}^{\pm}(z_1)\cdots
x^{\pm}_{i}(z_r)x^{\pm}_{j, l}x^{\pm}_{i}(z_{r+1})\cdots x^{\pm}_{i}(z_m)=0
\end{equation}
for any $l\in \mathbb Z$.
These identities are all the $w$-coefficients of (\ref{E:R9}).
\end{proof}

We now use the quantum vertex operator calculus \cite{J1} to prove
the Serre relations (\ref{E:R9}).

Repeatedly using (\ref{OPE2}) and Wick's theorem we have for
$m=1-(\a_i|\a_j)\geq 1$ (we only consider the $+$-case):
\begin{multline}\label{Ser1}
\sum_{r=0}^{m}\sum_{\sigma\in S_m}\sigma . \begin{bmatrix} m\\r\end{bmatrix}
X_i^+(z_1)\cdots X^+_i(z_r)X_j^+(w)X_i^+(z_{r+1})\cdots X_i(z_m)\\
=\sum_{r=0}^m\sum_{\sigma\in S_m}\sigma. \begin{bmatrix} m\\r\end{bmatrix}
 :X_i^+(z_1)\cdots X^+_i(z_r)X_j^+(w)X_i^+(z_{r+1})\cdots X_i(z_m):\\
\prod_{i<j}(z_i-q^{-1}z_j)_{q^2}^2\prod_{k=1}^{r}(z_k-q^{-1}w)_{q^2}^{1-m}
\prod_{k=r+1}^m(w-q^{-1}z_k)_{q^2}^{1-m}=0
\end{multline}
where the symmetric group $S_m$ acts on the indices of $z_i$'s:
$\sigma(z_i)=z_{\sigma(i)}$.

Pulling out the normal product $:X_j^+(x)X_i^+(z_1)\cdots X_i^+(z_m):$
and the anti-symmetric factor (under $S_m$ action):
$$\prod_{i<j}(z_i-z_j)\prod_{i=1}^m(w-z_i)_{q^2}^{2-m}\prod_{i=1}
(w-q^{1-m}z_i)^{-1}(z_i-q^{1-n}w)$$
and using
(\ref{OPE1}) we find that the Serre relation
is equivalent to the following identity:
\begin{multline} \label{Ser2}
\sum_{\sigma\in S_m}\sum_{r=0}^m(-1)^{l(\sigma)}\sigma.
\begin{bmatrix} m\\r\end{bmatrix}(w-q^{1-m}z_1)\cdots (w-q^{1-m}z_r)\\
\cdot (z_{r+1}-q^{1-m}w)\cdots (z_m-q^{1-m}w)\prod_{i<j}(z_i-q^{-2}z_j)=0
\end{multline}

Lemma (\ref{T:4}) says that we only need to show this identity for a special
component under the condition that there exists a root $\a_k$ joining with
$\a_j$ but not with $\a_i$. The constant coefficient of $w$ in (\ref{Ser2}) is
\begin{align*}
&\sum_{\sigma\in S_m}\sum_{r=0}(-1)^{l(\sigma)}\sigma.
\begin{bmatrix} m\\r\end{bmatrix}(-q^{1-m})^{r}z_1\cdots z_m
\prod_{i<j}(z_i-q^{-2}z_j) \\
=&z_1\cdots z_m(q^{-2m+2}; q^2)_m\sum_{\sigma\in S_m}(-1)^{l(\sigma)}
\sigma.\prod_{i<j}(z_i-q^{-2}z_j) \\
=& 0
\end{align*}
where we used the $q$-binomial theorem (\ref{binom1}).

Thus we have shown that our construction satisfies the Serre relation
under the special condition on $\a_i$ and $\a_j$, which in turn shows that
the combinatorial indentity (\ref{Ser2}) holds under this condition of
$\a_i$ and $\a_j$.  Note that the combinatorial identity
depends only on $(\a_i|\a_j)$, hence it is true in general.
Then we use identity (\ref{Ser2})
 to remove the restriction on $\a_i$ and $\a_j$.
Hence the Serre relation is proved in general, and our construction
indeed gives a representation of the $q$-affined Kac-Moody algebra, and
in turn a representation of the quantum Kac-Moody algebra.

\section{Combinatorial indentity}

In this paper we give a realization of the quantum Kac-Moody algebra
associated to the symmetric generalized Cartan matrix. As a by-product we
obtain the combinatorial identity (\ref{Ser2}), which seems to be
new. The indentity can be rewritten into the following form:

\begin{theorem} For any natural number $m$ we have:
\begin{multline} \label{Ser3}
\sum_{\sigma\in S_m}\sum_{r=0}^m\sigma.
\begin{bmatrix} m\\r\end{bmatrix}(w-q^{m-1}z_1)\cdots (w-q^{m-1}z_r)\\
\cdot (z_{r+1}-q^{m-1}w)\cdots (z_m-q^{m-1}w)\prod_{i<j}
\frac{z_i-q^{2}z_j}{z_i-z_j}=0
\end{multline}
\end{theorem}

Note that this gives an entirely non-trivial relation for Hall-Little\-wood
polynomials
defined in \cite{M}:

\begin{equation*}
P_{\lambda}(z; q^2)=\frac1{v_{\lambda}(q)}\sum_{\sigma\in S_m}\sigma.(
z_1^{\lambda_1}\cdots z_m^{\lambda_m}\prod_{i<j}
\frac{z_i-q^{2}z_j}{z_i-z_j}),
\end{equation*}
where $\lambda$ is a partition of $m$, and $v_{\lambda}(q)=
q^{\sum\binom{\lambda_i}{2}}\prod_i[\lambda_i]!$
. The collection of $\{P_{\lambda}|\lambda\
\mbox{partitions of}\ $m$\}$ forms a $\mathbb Z$-basis for the ring
of symmetric functions of degree $m$. With some effort the notion of
$P_{\lambda}$ can be extended to $m$-tuples.
Our identity gives a relation
for all $P_{\lambda}$ when $\lambda$ are not necessary partitions.

 Although we proved the identity (\ref{Ser3}) by our representation of
the quantum Kac-Moody algebra, it would be interesting to prove this
 combinatorial identity independently from quantum groups, which will
be treated elsewhere.

The method used in this paper can be generalized to higher level cases
using the quantum $Z$-algebra method (cf. \cite{J4}).

\end{document}